\documentclass{article}
\usepackage{amssymb}
\usepackage{graphicx}
\usepackage{amsmath}
\usepackage{array}

\newcommand{\A}{{\mathcal A}}
\newcommand{\B}{{\mathcal B}}
\newcommand{\C}{{\mathcal C}}
\newcommand{\D}{{\mathcal D}}
\newcommand{\E}{{\mathcal E}}
\newcommand{\I}{{\mathcal I}}
\renewcommand{\L}{{\mathcal L}}
\renewcommand{\S}{{\mathcal S}}
\newcommand{\T}{{\mathcal T}}

\newcommand{\Chi}{\overline \chi}

\newcommand{\kk}{\Bbbk}
\newcommand{\RR}{\mathbb{R}}
\newcommand{\FF}{\mathbb{F}}
\newcommand{\NN}{\mathbb{N}}
\newcommand{\FFr}{{\mathbb{F}_q^{\,r}}}
\newcommand{\FFn}{{\mathbb{F}_q^{\,n}}}
\newcommand{\ZZ}{\mathbb{Z}}

\newtheorem{theorem}{Theorem}[section]
\newtheorem{defin}[theorem]{Definition}
\newtheorem{proposition}[theorem]{Proposition}

\numberwithin{equation}{section}

\renewcommand{\dim}{\mathop{\rm dim}\nolimits}

\renewcommand{\Im}{\mathop{\rm Im}\nolimits}

\begin{document}

\pagestyle{plain}

\title{Computing the Tutte polynomial \\ of a hyperplane
arrangement}
\author{Federico Ardila}
\date{}
%

\maketitle

\begin{abstract}
We define and study the Tutte polynomial of a hyperplane
arrangement. We introduce a method for computing it by solving an
enumerative problem in a finite field. For specific arrangements,
the computation of Tutte polynomials is then reduced to certain
related enumerative questions. As a consequence, we obtain new
formulas for the generating functions enumerating alternating
trees, labelled trees, semiorders and Dyck paths.
\end{abstract}

\section{Introduction.}\label{sec:introarr}

Much work has been devoted in recent years to studying hyperplane
arrangements and, in particular, their characteristic polynomials.
The polynomial $\chi_{\A}(q)$ is a very powerful invariant of the
arrangement $\A$; it arises very naturally in many different
contexts. Two of the many beautiful results about the
characteristic polynomial of an arrangement are the following.

\begin{theorem}(Zaslavsky, \cite{Za75})
Let $\A$ be a hyperplane arrangement in $\RR^n$. The number of
regions into which $\A$ dissects $\RR^n$ is equal to $(-1)^n
\chi_{\A}(-1)$. The number of regions which are relatively
bounded is equal to $(-1)^n \chi_{\A}(1)$.
\end{theorem}

\begin{theorem}(Orlik-Solomon, \cite{Or80})
Let $\A$ be a hyperplane arrangement in $\mathbb{C}^n$, and let
$M_{\A} = \mathbb{C}^n - \bigcup_{H \in \A} H$ be its complement.
Then the Poincar\'e polynomial of the cohomology ring of $M_{\A}$
is given by:
$$
\sum_{k \geq 0} \mathrm{rank} \, H^k(M_{\A}, \ZZ) \, q^k = (-q)^n
\chi_{\A}(-1/q).
$$
\end{theorem}

Several authors have worked on computing the characteristic
polynomials of specific hyperplane arrangements. This work has
led to some very nice enumerative results; see for example
\cite{At96}, \cite{Po00}.

\medskip

It is somewhat surprising that nothing has been said about the
Tutte polynomial of a hyperplane arrangement. Graphs and matroids
have a Tutte polynomial associated with them, which generalizes
the characteristic polynomial, and arises very naturally in
numerous enumerative problems in both areas. Many interesting
invariants of graphs and matroids can be computed immediately from
this polynomial.


The present paper, in conjunction with \cite{semimatroids}, aims
to define and investigate the Tutte polynomial of a hyperplane
arrangement. This paper is devoted to purely enumerative
questions. We are particularly interested in computing the Tutte
polynomials of specific arrangements. We address the
matroid-theoretic aspects of this investigation in
\cite{semimatroids}.


The paper is organized as follows. In Section
\ref{sec:arrangements} we introduce the basic notions of
hyperplane arrangements that we will need. In Section
\ref{sec:compTutte} we define the Tutte polynomial of a hyperplane
arrangement, and we present a finite field method for computing
it.  This is done in terms of the coboundary polynomial, a simple
transformation of the Tutte polynomial. We recover several known
results about the characteristic and Tutte polynomials of graphs
and representable matroids, and derive other consequences of this
method. Finally, in Section \ref{sec:comp}, we compute the Tutte
polynomials of several families of arrangements. In particular,
for deformations of the braid arrangement, we relate the
computation of Tutte polynomials to the enumeration of classical
combinatorial objects. As a consequence, we obtain several purely
enumerative results about objects such as labeled trees, Dyck
paths, alternating trees and semiorders.

\section{Hyperplane arrangements.}\label{sec:arrangements}

In this section we recall some of the basic concepts of hyperplane
arrangements. For a more thorough introduction, we refer the
reader to \cite{Or92} or \cite{St04}.

Given a field $\kk$ and a positive integer $n$, an \emph{affine
hyperplane} in $\kk^n$ is an $(n-1)$-dimensional affine subspace
of $\kk^n$. If we fix a system of coordinates $x_1, \ldots, x_n$
on $\kk^n$, a hyperplane can be seen as the set of points that
satisfy a certain equation $c_1x_1 + \cdots + c_nx_n = c$, where
$c_1, \ldots, c_n, c \in \kk$ and not all $c_i$'s are equal to
$0$. A \emph{hyperplane arrangement} $\A$ in $\kk^n$ is a finite
set of affine hyperplanes of $\kk^n$. We will refer to hyperplane
arrangements simply as \emph{arrangements}. We will assume for
simplicity that $\kk = \RR$ unless explicitly stated, although
most of our results extend immediately to any field of
characteristic zero.

We will say that an arrangement $\A$ is \emph{central} if the
hyperplanes in $\A$ have a non-empty intersection.\footnote{
Sometimes we will call an arrangement \emph{affine} to emphasize
that it does not need to be central.} Similarly, we will say that
a subset (or \emph{subarrangement}) $\B \subseteq \A$ of
hyperplanes is \emph{central} if the hyperplanes in $\B$ have a
non-empty intersection.

The \emph{rank function} $r_{\A}$ is defined for each central
subset $\B$ by $r_{\A}(\B) = n - \dim \cap \B$. This function can
be extended to a function $r_{\A}:2^{\A} \rightarrow \NN$, by
defining the rank of a non-central subset $\B$ to be the largest
rank of a central subset of $\B$. The \emph{rank} of $\A$ is
$r_{\A}(\A)$, and it is denoted $r_{\A}$.

Alternatively, if the hyperplane $H$ has defining equation $c_1x_1
+ \cdots +c_nx_n = c$, associate its normal vector $v = (c_1,
\ldots, c_n)$ to it. Then define $r_{\A}(\{H_1, \ldots ,H_k\})$ to
be the dimension of the span of the corresponding vectors $v_1,
\ldots, v_k$ in $\RR^n$. It is easy to see that these two
definitions of the rank function agree.  In particular, this means
that the resulting function $r_{\A}:2^{\A} \rightarrow \NN$ is the
rank function of a matroid. We will usually omit the subscripts
when the underlying arrangement is clear, and simply write $r(\B)$
and $r$ for $r_{\A}(\B)$ and $r_{\A}$,
respectively.\footnote{There is another natural way to extend
$r_{\A}$ to the rank function of a matroid; for more information,
see \cite{semimatroids}.}

The rank function gives us natural definitions of the usual
concepts of matroid theory, such as independent sets, bases,
closed sets, and circuits, in the context of hyperplane
arrangements. All of this is done more naturally in the broader
context of semimatroids in \cite{semimatroids}.


To each hyperplane arrangement $\A$ we assign a partially ordered
set, called the \emph{intersection poset} of $\A$ and denoted
$L_{\A}$. It consists of the non-empty intersections $H_{i_1} \cap
\cdots \cap H_{i_k}$, ordered by reverse inclusion. This poset is
graded, with rank function $r(H_{i_1} \cap \cdots \cap H_{i_k}) =
r_{\A}(\{H_{i_1}, \ldots, H_{i_k}\})$, and a unique minimal
element $\hat{0} = \RR^n$. We will sometimes call two arrangements
$\A_1$ and $\A_2$ \emph{isomorphic}, and write $\A_1 \cong \A_2$,
if $L_{\A_1} \cong L_{\A_2}$.

The \emph{characteristic polynomial} of $\A$ is
$$
\chi_{\A}(q) = \sum_{x \in L_{\A}} \mu(\hat{0}, x) q^{n - r(x)}.
$$
where $\mu$ denotes the M\"{o}bius function \cite[Section
3.7]{St86} of $L_{\A}$.

Let $\A$ be an arrangement and let $H$ be a hyperplane in $\A$.
The arrangement $\A - \{H\}$ (or simply $\A-H$), obtained by
removing $H$ from the arrangement, is called the \emph{deletion of
$H$ in $\A$}. It is an arrangement in $\RR^n$. The arrangement
$\A/H = \{H' \cap H \, |\, H'\in \A-H, H' \cap H \neq
\emptyset\}$, consisting of the intersections of the other
hyperplanes with $H$, is called the \emph{contraction of $H$ in
$\A$}. It is an arrangement in $H$.

Notice, however, that some technical difficulties can arise. In a
hyperplane arrangement $\A$, contracting a hyperplane $H$ may give
us repeated hyperplanes $H_1$ and $H_2$ in the arrangement $\A/H$.
Now say we want to contract $H_1$ in $\A/H$. In passing to the
contraction $(\A/H)/H_1$, the hyperplane $H_2$ of $\A/H$ becomes
the "hyperplane" $H_2 \cap H_1 = H_1$ in the ``arrangement"
$(\A/H)/H_1$. This is not a hyperplane in $H_1$, though.

Therefore, the class of hyperplane arrangements, as we defined it,
is not closed under deletion and contraction. This is problematic
when we want to mirror matroid-theoretic results in this context.
There is an artificial solution to this problem: we can consider
multisets $\{H_1, \ldots, H_k\}$ of subspaces of vector spaces
$V$, where each $H_i$ has dimension $\dim V-1$ or $\dim V$. In
other words, we allow repeated hyperplanes, and we allow the full
space $V$ to be regarded as a ``hyperplane", mirroring a loop of a
matroid. This class of objects \emph{is} closed under deletion and
contraction, but it is somewhat awkward to work with. A better
solution is to think of arrangements as members of the class of
semimatroids; a class that \emph{is} closed under deletion and
contraction, and is more natural matroid-theoretically. We develop
this point of view in \cite{semimatroids}. However, such issues
will be irrelevant in this paper, which focuses on purely
enumerative aspects of arrangements.

\section{Computing the Tutte polynomial.}\label{sec:compTutte}

In \cite{At96}, Athanasiadis introduced a powerful method for
computing the characteristic polynomial of a subspace arrangement,
based on ideas of Crapo and Rota \cite{Cr70}. He reduced the
computation of characteristic polynomials to an enumeration
problem in a vector space over a finite field. He used this method
to compute explicitly the characteristic polynomial of several
families of hyperplane arrangements, obtaining very nice
enumerative results. As should be expected, this method only works
when the equations defining the hyperplanes of the arrangement
have integer (or rational) coefficients. Such an arrangement will
be called a \emph{$\ZZ-$arrangement}.

In \cite{Re99}, Reiner asked whether it is possible to use
\cite[Corollary 3]{Re99} to compute explicitly the Tutte
polynomials of some non-trivial families of representable
matroids. Compared to all the work that has been done on computing
characteristic polynomials explicitly, virtually nothing is known
about computing Tutte polynomials.

In this section, we introduce a new method for computing Tutte
polynomials of hyperplane arrangements. Our approach does not use
Reiner's result; it is closer to Athanasiadis's method. In fact,
Athanasiadis's result \cite[Theorem 2.2]{At96} can be obtained as
a special case of the main result of this section, Theorem
\ref{th:finfield}, by setting $t=0$.

After proving Theorem \ref{th:finfield}, we will present some of
its consequences. We will then use it in Section \ref{sec:comp}
to compute explicitly the Tutte polynomials of several families
of arrangements.

\subsection{The Tutte and coboundary polynomials.}\label{sec:cob}

\begin{defin}\label{def:Tuttehyp}
The \emph{Tutte polynomial} of a hyperplane arrangement $\A$ is
\begin{equation}\label{eq:Tuttehyp}
T_{\A}(q,t) = \sum_{\stackrel{\B \subseteq \A}{\mathrm{central}}}
(q-1)^{r-r(\B)}(t-1)^{|\B|-r(\B)},
\end{equation}
where the sum is over all central subsets $\B \subseteq \A$.
\end{defin}

It will be useful for us to consider a simple transformation of
the Tutte polynomial, first considered by Crapo \cite{Cr69} in
the context of matroids.

\begin{defin}
The \emph{coboundary polynomial} $\Chi_{\A}(q,t)$ of an
arrangement $\A$ is
\begin{equation}\label{eq:cob}
\Chi_{\A}(q,t) = \sum_{\stackrel{\B \subseteq
\A}{\mathrm{central}}} q^{r-r(\B)} (t-1)^{|\B|}.
\end{equation}
\end{defin}

It is easy to check that
$$
\Chi_{\A}(q,t) = (t-1)^r \, T_{\A}\left(\frac{q+t-1}{t-1},
t\right)
$$
and
$$
T_{\A}(x,y) = \frac1{(y-1)^r} \,
\Chi_{\A}\left((x-1)(y-1),y\right).
$$

Therefore, computing the coboundary polynomial of an arrangement
is essentially equivalent to computing its Tutte polynomial. Our
results can be presented more elegantly in terms of the coboundary
polynomial.

\subsection{The finite field method.}\label{sec:finfieldform}

Let $\A$ be a $\ZZ-$arrangement in $\RR^n$, and let $q$ be a prime
power. The arrangement $\A$ induces an arrangement $\A_q$ in the
vector space $\FFn$. If we consider the equations defining the
hyperplanes of $\A$, and regard them as equations over $\FF_q$,
they define the hyperplanes of $\A_q$.

Say that $\A$ \emph{reduces correctly} over $\FF_q$ if the
arrangements $\A$ and $\A_q$ are isomorphic. This does not always
happen; sometimes the hyperplanes of $\A$ do not even become
hyperplanes in $\A_q$. For example, the hyperplane $2x+2y=1$ in
$\RR^2$ becomes the empty ``hyperplane" $0=1$ in $\FF_2^{\,2}$.
Sometimes independence is not preserved. For example, the
independent hyperplanes $2x+y=0$ and $y=0$ in $\RR^2$ become the
same hyperplane in $\FF_2^{\,2}.$

However, if $q$ is a power of a large enough prime, $\A$ will
reduce correctly over $\FF_q$. To have $\A \cong \A_q$, we need
central and independent subarrangements to be preserved. Cramer's
rule lets us rephrase these conditions, in terms of certain
determinants (formed by the coefficients of the hyperplanes in
$\A$)  being zero or non-zero. If we let $q$ be a power of a prime
which is larger than all these determinants, we will guarantee
that $\A$ reduces correctly over $\FF_q$.

\begin{theorem}\label{th:finfield}
Let $\mathcal A$ be a $\ZZ-$arrangement in $\RR^n$. Let $q$ be a
power of a large enough prime, and let ${\A}_q$ be the induced
arrangement in $\FFn$. Then
\begin{equation}\label{eq:finfield}
q^{n-r} \, \Chi_{\A}(q,t) \, = \sum_{p \, \in \, \FFn} t^{h(p)}
\end{equation}
where $h(p)$ denotes the number of hyperplanes of ${\A}_q$ that
$p$ lies on.
\end{theorem}

\noindent \emph{Proof.} Let $q$ be a power of a large enough
prime, so that $\A$ reduces correctly over $\FF_q$. For each $\B
\subseteq \A$, let $\B_q$ be the subarrangement of $\A_q$ induced
by it. For each $p \in \FFn$, let $H(p)$ be the set of hyperplanes
of $\A_q$ that $p$ lies on. From (\ref{eq:cob}) we have

\begin{eqnarray*}
q^{n-r}\,\Chi_{\A}(q,t) &=& \sum_{\stackrel{\B \subseteq
\A}{\mathrm{central}}} q^{n-r(\B)} (t-1)^{|\B|} \\
& = & \sum_{\stackrel{\B \subseteq \A}{\mathrm{central}}}
q^{\dim \cap \B} (t-1)^{|\B|} \\
& = & \sum_{\stackrel{\B_q \subseteq \A_q}{\mathrm{central}}}
|\cap \B_q| \, (t-1)^{|\B_q|} \\
& = & \sum_{\stackrel{\B_q \subseteq \A_q}{\mathrm{central}}}
\sum_{p \, \in \, \cap \B_q} \, (t-1)^{|\B_q|} \\
& = & \sum_{p \, \in \, \FFn} \sum_{\B_q \subseteq H(p)}
(t-1)^{|\B_q|} \\
& = & \sum_{p \, \in \, \FFn} (1+(t-1))^{h(p)},
\end{eqnarray*}
as desired. $\Box$

\medskip

In principle, Theorem \ref{th:finfield} only computes
$\Chi_{\A}(q,t)$ when $q$ is a power of a large enough prime. In
practice, however, when we compute the right-hand side of
(\ref{eq:finfield}) for large prime powers $q$, we will get a
polynomial function in $q$ and $t$. Since the left-hand side is
also a polynomial, these two polynomials must be equal.

Theorem \ref{th:finfield} reduces the computation of coboundary
polynomials (and hence Tutte polynomials) to enumerating points
in the finite vector space $\FFn$, according to a certain
statistic. This method can be extremely useful when the
hyperplanes of the arrangement are defined by simple equations.
We will illustrate this in section \ref{sec:comp}.

We remark that Theorem \ref{th:finfield} was also obtained
independently by Welsh and Whittle \cite{We02}.

\subsection{Special cases and applications.}\label{sec:apps}

Now we present some known facts and some new results which follow
from the finite field method. We start with two classical theorems
which are special cases of Theorem \ref{th:finfield}.

\subsubsection{Colorings of graphs}

A graph $G$ has a matroid $M_G$ associated to it, called the
\emph{cycle matroid} of $G$. Its Tutte polynomial is equal to the
(graph-theoretic) Tutte polynomial of $G$.

From the point of view of arrangements, the construction is the
following. Given a graph $G$ on $[n]$, we associate to it an
arrangement $\A_G$ in $\RR^n$. It consists of the hyperplanes
$x_i = x_j$, for all $1 \leq i < j \leq n$ such that $ij$ is an
edge in the graph $G$. Then we have that $T_G(q,t) =
T_{\A_G}(q,t)$. We can define the coboundary polynomial for a
graph like we did for arrangements, and then $\Chi_G(q,t) =
\Chi_{\A_G}(q,t)$ also.

We shall now interpret Theorem \ref{th:finfield} in this
framework. It is easy to see that the rank of $G$ is equal to
$n-c$, where $c$ is the number of connected components of $G$.
Therefore the left-hand side of (\ref{eq:finfield}) is $q^c\,
\Chi_{G}(q,t)$ in this case.

To interpret the right-hand side, notice that each point $p \in
\FFn$ corresponds to a $q$-coloring of the vertices of $G$. The
point $p=(p_1, \ldots, p_n)$ will correspond to the coloring
$\kappa_p$ of $G$ which assigns color $p_i$ to vertex $i$. A
hyperplane $x_i=x_j$ contains $p$ when $p_i=p_j$. This happens
precisely when edge $ij$ is \emph{monochromatic} in $\kappa_p$;
that is, when its two ends have the same color. Therefore,
applying Theorem \ref{th:finfield} to the arrangement $\A_G$, we
recover the following known result:

\begin{theorem}(\cite[Proposition 6.3.26]{Br92})
Let $G$ be a graph with $c$ connected components. Then
$$
q^c \, \Chi_G(q,t) = \sum_{\stackrel{q-\mathrm{colorings}}{\kappa
\, \mathrm{of} \, G}} t^{ {\mathop{\rm mono}\nolimits}(\kappa) },
$$
where ${\mathop{\rm mono}\nolimits}(\kappa)$ is the number of
monochromatic edges in $\kappa$.
\end{theorem}

\subsubsection{Linear codes}

Given positive integers $n \geq r$, an \emph{$[n,r]$ linear code}
$C$ over $\FF_q$ is an $r$-dimensional subspace of $\FFn$. A
\emph{generator matrix} for $C$ is an $r \times n$ matrix $U$
over $\FF_q$, the rows of which form a basis for $C$. It is not
difficult to see that the isomorphism class of the matroid on the
columns of $U$ depends only on $C$. We shall denote the
corresponding matroid $M_C$.

The elements of $C$ are called \emph{codewords}. The \emph{weight}
$w(v)$ of a codeword is the cardinality of its support; that is,
the number of non-zero coordinates of $v$. The \emph{codeweight
polynomial} of $C$ is
\begin{equation}\label{eq:codeweight}
A(C,q,t) = \sum_{v \, \in \, C} t^{w(v)}.
\end{equation}

The translation of Theorem \ref{th:finfield} to this setting is
the following.

\begin{theorem}(Greene, \cite{Gr76})
For any linear code $C$ over $\FF_q$,
$$
A(C, q, t) = t^n \, \Chi_{M_C}\left(q, \frac1t\right).
$$
\end{theorem}

\noindent \emph{Proof.} Let $\A_C$ be the central arrangement
corresponding to the columns of $U$. (We can call it $\A_C$
because, as stated above, its isomorphism class depends only on
$C$.) This is a rank $r$ arrangement in $\FFr$ such that
$\Chi_{M_C}(q, \frac1t) = \Chi_{\A_C}(q, \frac1t)$. Comparing
(\ref{eq:codeweight}) with Theorem \ref{th:finfield}, it remains
to prove that
$$
\sum_{v \, \in \, C} t^{w(v)} = \sum_{p \, \in \, \FFr}
t^{n-h(p)}.
$$

To do this, consider the bijection $\phi:\FFr \rightarrow C$
determined by right multiplication by $U$. If $u_1, \ldots, u_r$
are the row vectors of $U$, then $\phi$ sends $p=(p_1, \ldots,
p_r) \in \FFr$ to the codeword $v_p = p_1u_1 + \cdots + p_ru_r \in
C$. For $1 \leq i \leq n$, $p$ lies on the hyperplane determined
by the $i-th$ column of $U$ if and only if the $i$-th coordinate
of $v_p$ is equal to zero. Therefore $h(p) = n - w(v_p)$. This
completes the proof. $\Box$

\subsubsection{Deletion-contraction}\label{sec:d-c rev}

The point of view of Theorem \ref{th:finfield} can be used to give
a nice enumerative proof of the deletion-contraction formula for
the Tutte polynomial of an arrangement. Once again, this formula
is better understood in the context of semimatroids, as shown in
\cite{semimatroids}. For the moment, leaving matroid-theoretical
issues aside, we only wish to present a special case of it as a
nice application.

\begin{proposition}
Let $\A$ be a hyperplane arrangement, and let $H$ be a hyperplane
in $\A$ such that $r_{\A}(\A-H) = r_{\A}$. Then $T_{\A}(q,t) =
T_{\A-H}(q,t) + T_{\A/H}(q,t)$.
\end{proposition}

\medskip

\noindent \emph{Proof.} Because there will be several
arrangements involved, let $h(\B,p)$ denote the number of
hyperplanes in $\B_q$ that $p$ lies on. Then

\begin{eqnarray*}
q^{n-r}\,\Chi_{\A}(q,t) &=& \sum_{p \, \in \, \FFn} t^{h(\A, \, p)} \\
&=& \sum_{p \, \in \, \FFn-H} t^{h(\A, \, p)} + \sum_{p \, \in
\, H} t^{h(\A, \, p)} \\
&=& \sum_{p \, \in \, \FFn-H} t^{h(\A-H,\,p)} + \sum_{p \, \in \,
H} t^{h(\A-H,\,p)+1} \\
&=& \sum_{p \, \in \, \FFn} t^{h(\A-H,\,p)} +(t-1)\, \sum_{p \,
\in \, H}
t^{h(\A-H,\,p)} \\
&=& q^{n-r}\,\Chi_{\A-H}(q,t) + (t-1)q^{(n-1)-(r-1)}
\,\Chi_{\A/H}(q,t).
\end{eqnarray*}
We conclude that $\Chi_{\A}(q,t) = \Chi_{\A-H}(q,t) + (t-1) \,
\Chi_{\A/H}(q,t)$, which is equivalent to the
deletion-contraction formula for Tutte polynomials. $\Box$

\subsubsection{A probabilistic interpretation}

\begin{theorem}
Let $\A$ be an arrangement and let $0 \leq t \leq 1$ be a real
number. Let $\B$ be a random subarrangement of $\A$, obtained by
independently removing each hyperplane from $\A$ with probability
$t$. Then the expected characteristic polynomial $\chi_{\B}(q)$
of $\B$ is $q^{n-r}\Chi_{\A}(q,t)$.
\end{theorem}

\noindent \emph{Proof.} We have
\begin{eqnarray*}
E[\chi_{\B}(q)] &=& \sum_{\C \subseteq \A} P[\B = \C] \,
\chi_{\C}(q)
\\
&=& \sum_{\C \subseteq \A} P[\B = \C] \, |\, \FFn - \,\, \cup \C_q| \\
&=& \sum_{p \, \in \, \FFn} \sum_{\stackrel{\C \subseteq \A}{p \,
\notin \, \cup \, \C_q}} P[\B = \C],
\end{eqnarray*}
where in the second step we have used Athanasiadis's result
\cite{At96}; that is, the case $t=0$ of Theorem \ref{th:finfield}.

Recall that $H(p)$ denotes the set of hyperplanes in $\A_q$
containing $p$. Then
\begin{eqnarray*}
E[\chi_{\B}(q)] &=& \sum_{p \, \in \, \FFn} P[\B_q \cap H(p) = \emptyset] \\
&=& \sum_{p \, \in \, \FFn} t^{h(p)},
\end{eqnarray*}
which is precisely what we wanted to show. $\Box$

\subsubsection{A M\"{o}bius formula}

\begin{theorem}
For an arrangement $\A$ and an affine subspace $x$ in the
intersection poset $L_{\A}$, let $h(x)$ be the number of
hyperplanes of $\A$ containing $x$. Then
$$
\Chi_{\A}(q,t) = \sum_{x \leq y \,\, \mathrm{in} \, L_{\A}}
\mu(x,y)\, q^{r-r(y)} t^{h(x)} .
$$
\end{theorem}

\noindent \emph{Proof.} Consider the arrangement $\A$ restricted
to $\FFn$, where $q$ is a power of a large enough prime, so that
$\A$ reduces correctly over $\FF_q$. Given $x \in L_{\A_q}$, let
$P(x)$ be the set of points in $\FFn$ which are contained in $x$,
and are not contained in any $y$ such that $y > x$ in $L_{\A_q}$.
Then the set $x$ is partitioned by the sets $P(y)$ for $y \geq
x$, so we have
$$
q^{\dim x} = |x| = \sum_{y \geq x} |P(y)|.
$$
By the M\"{o}bius inversion formula \cite[Proposition 3.7.1]{St86}
we have
$$
|P(x)| = \sum_{y \geq x} \mu(x,y) \, q^{\dim y}.
$$
Now, from Theorem \ref{th:finfield} we know that
\begin{eqnarray*}
q^{n-r}\Chi_{\A}(q,t) &=& \sum_{x \in L_{\A}} \sum_{p \in P(x)}
t^{h(p)} = \sum_{x \in L_{\A}} |P(x)| \, t^{h(x)} \\
&=& \sum_{x \leq y \,\, \mathrm{in} \, L_{\A}} \mu(x,y) \,
q^{n-r(y)} t^{h(x)},
\end{eqnarray*}
as desired. $\Box$

\section{Computing coboundary polynomials.}\label{sec:comp}

In this section we use Theorem \ref{th:finfield} to compute the
coboundary polynomials of several families of arrangements. As
remarked at the beginning of Section \ref{sec:cob}, this is
essentially the same as computing their Tutte polynomials.

\subsection{Coxeter arrangements.}\label{sec:cox}

To illustrate how our finite field method works, we start by
presenting some simple examples.

Let $\Phi$ be an irreducible crystallographic root system in
$\RR^n$, and let $W$ be its associated Weyl group. The Coxeter
arrangement of type $W$ consists of the hyperplanes $(\alpha,x) =
0$ for each $\alpha \in \Phi^+$, with the standard inner product
on $\RR^n$. See \cite{Hu90} for an introduction to root systems
and Weyl groups, and \cite[Chapter 6]{Or92} or \cite[Section
2.3]{Bj93} for more information on Coxeter arrangements.

In this section we compute the coboundary polynomials of the
Coxeter arrangements of type $A_n$,
$B_n$ and $D_n$. (The arrangement of type $C_n$ is the
same as the arrangement of type $B_n$.) The best way
to state our results is to compute the exponential generating
function for the coboundary polynomials of each family.

The following three theorems have never been stated explicitly in
the literature in this form. Theorem \ref{th:An} is equivalent to
a result of Tutte \cite{Tu53}, who computed the Tutte polynomial
of the complete graph. It is also an immediate consequence of a
more general theorem of Stanley \cite[(15)]{St98b}. Theorems
\ref{th:Bn} and \ref{th:Dn} are implicit in the work of Zaslavsky
\cite{Za95}.

\begin{theorem}\label{th:An}
Let $\A_n$ be the Coxeter arrangement of type $A_{n-1}$ in
$\RR^n$, consisting of the hyperplanes $x_i=x_j$ for $1 \leq i <
j \leq n$.\footnote{This arrangement is also known as the
\emph{braid arrangement.}} We have
$$
1 + q \sum_{n \geq 1} \Chi_{\A_n}(q,t)\frac{x^n}{n!} =
\left( \sum_{n \geq 0} t^{n \choose 2} \frac{x^n}{n!} \right)^q.
$$
\end{theorem}

\noindent \emph{Proof.} For $n \geq 1$ we have that
$$
q  \, \Chi_{\A_n}(q,t) = \sum_{p \, \in \, \FFn} t^{h(p)}.
$$
for all powers of a large enough prime $q$, according to Theorem
\ref{th:finfield}. For each $p \in \FFn$, if we let $A_k = \{i \in
[n] \, | \, p_i = k\}$ for $0 \leq k \leq q-1$, then $h(p) =
{|A_0| \choose 2} + \cdots + {|A_{q-1}| \choose 2}.$ Thus
$$
q \, \Chi_{\A_n}(q,t) = \sum_{A_0 \cup \, \cdots \, \cup A_{q-1} =
[n]} t^{{|A_0| \choose 2}+ \, \cdots \, + {|A_{q-1}| \choose 2}}
$$
where the sum is over all weak ordered $q$-partitions of $[n]$.
The compositional formula for exponential generating functions
\cite[Proposition 5.1.3]{St99}, \cite{Be97} implies the desired
result. $\Box$

\medskip

\begin{theorem}\label{th:Bn}
Let $\B_n$ be the Coxeter arrangement of type $B_n$ in $\RR^n$,
consisting of the hyperplanes $x_i = x_j$ and $x_i + x_j=0$ for
$1 \leq i < j \leq n$, and the hyperplanes $x_i = 0$ for $1 \leq
i \leq n$. We have
$$
\sum_{n \geq 0} \Chi_{\B_n}(q,t) \frac{x^n}{n!} = \left(\sum_{n
\geq 0} 2^n t^{n \choose 2} \frac{x^n}{n!} \right)^{\frac{q-1}2}
\left(\sum_{n \geq 0} t^{\, n^2} \frac{x^n}{n!} \right).
$$
\end{theorem}

\noindent \emph{Proof.} Let $q$ be a power of a large enough
prime, and  let $s = \frac{q-1}2.$ Now for each $p \in \FFn$, if
we let $B_k = \{i \in [n] \, | \, p_i = k \, \mathrm{\, or \,} \,
p_i = q-k\}$ for $0 \leq k \leq s$, we have that $h(p) =
|B_0|^{\,2} + {|B_1| \choose 2} + \cdots + {|B_s| \choose 2}$.
Also, given a weak ordered partition $(B_0, \ldots, B_s)$ of
$[n]$, there are $2^{|B_1| + \cdots + |B_s|}$ points of $p$ which
correspond to it: for each $i \in B_k$ with $k \neq 0$, we get to
choose whether $p_i$ is equal to $k$ or to $q-k$. Therefore
$$
q \, \Chi_{\B_n}(q,t) = \sum_{B_0 \cup \, \cdots \, \cup B_s =
[n]} t^{|B_0|^2} \left(2^{|B_1|}t^{|B_1| \choose 2}\right) \cdots
\left(2^{|B_s|}t^{|B_s| \choose 2}\right),
$$
and the compositional formula for exponential generating functions
implies Theorem \ref{th:Bn}. $\Box$ \medskip

\begin{theorem}\label{th:Dn}
Let $\D_n$ be the Coxeter arrangement of type $D_n$ in $\RR^n$,
consisting of the hyperplanes $x_i = x_j$ and $x_i + x_j=0$ for
$1 \leq i < j \leq n$. We have
$$
\sum_{n \geq 0} \Chi_{\D_n}(q,t) \frac{x^n}{n!} = \left(\sum_{n
\geq 0} 2^nt^{n \choose 2} \frac{x^n}{n!} \right)^{\frac{q-1}2}
\left(\sum_{n \geq 0} t^{\, n(n-1)} \, \frac{x^n}{n!} \right).
$$
\end{theorem}

\medskip

We omit the details of the proof of Theorem \ref{th:Dn}, which is
very similar to the proof of Theorem \ref{th:Bn}.

\medskip

Setting $t=0$ in Theorems \ref{th:An}, \ref{th:Bn} and \ref{th:Dn},
it is easy to recover the well-known formulas for the characteristic
polynomials of the above arrangements:
\begin{eqnarray*}
\chi_{\A_n}(q) & = & q(q-1)(q-2) \cdots (q-n+1), \\
\chi_{\B_n}(q) & = & (q-1)(q-3) \cdots (q-2n+1), \\
\chi_{\D_n}(q) & = & (q-1)(q-3) \cdots (q-2n+3)(q-n+1).
\end{eqnarray*}

\subsection{Two more examples.}

\begin{theorem}\label{th:An'}
Let ${\mathcal A}_n^{\#}$ be a generic deformation of the
arrangement $\A_n$, consisting of the hyperplanes $x_i - x_j =
a_{ij}$ $(1 \leq i < j \leq n)$, where the $a_{ij}$ are generic
real numbers \footnote{The $a_{ij}$ are ``generic" if no $n$ of
the hyperplanes have a non-empty intersection, and any non-empty
intersection of $k$ hyperplanes has rank $k$. This can be
achieved, for example, by requiring that the $a_{ij}$'s are
linearly independent over the rational numbers. Almost all choices
of the $a_{ij}$'s are generic.}. For $n \geq 1$,
$$
q \,\Chi_{\A_n^{\#}} (q,t) = \sum_F q^{n - e(F)}(t-1)^{e(F)}
$$
where the sum is over all forests $F$ on $[n]$, and $e(F)$
denotes the number of edges of $F$. Also,
$$
1 + q\,\sum_{n \geq 1} \Chi_{{\mathcal A}^{\#}_n}(q,t)
\frac{x^n}{n!} = \left(\sum_{n \geq 0} f(n) \frac{x^n(t-1)^n}{n!}
\right)^{\frac{q}{t-1}},
$$
where $f(n)$ is the number of forests on $[n]$.
\end{theorem}

\noindent \emph{Proof.} It is possible to prove Theorem
\ref{th:An'} using our finite field method, as we did in the
previous section. However, it will be easier to proceed directly
from (\ref{eq:cob}), the definition of the coboundary polynomial.

To each subarrangement $\B$ of $\A_n^{\#}$ we can assign a graph
$G_{\B}$ on the vertex set $[n]$, by letting edge $ij$ be in
$G_{\B}$ if and only if the hyperplane $x_i - x_j = a_{ij}$ is in
$\B$. Since the $a_{ij}$'s are generic, the subarrangement $\B$
is central if and only if the corresponding graph $G_{\B}$ is a
forest. For such a $\B$, it is clear that $|\B| = r(\B) =
e(G_{\B})$. Hence,
\begin{eqnarray*}
\Chi_{\A_n^{\#}}(q,t) &=& \sum_{\stackrel{\B \subseteq \A_n^{\#}}
{\mathrm{central}}} q^{r - r(\B)} (t-1)^{|\B|}\\
&=& \sum_F q^{(n-1)-e(F)}(t-1)^{e(F)},
\end{eqnarray*}
proving the first claim. Now let $c(F)=n-e(F)$ be the number of
connected components of $F$. We have
\begin{eqnarray*}
1 + q\,\sum_{n \geq 1} \Chi_{{\mathcal A}^{\#}_n}(q,t)
\frac{x^n}{n!} &=& \sum_{n \geq 0} \sum_{F \, \mathrm{on} \, [n]}
\left(\frac{q}{t-1}\right)^{c(F)} \frac{x^n(t-1)^n}{n!} \\
&=&\left(\sum_{n \geq 0} f(n) \frac{x^n(t-1)^n}{n!}
\right)^{\frac{q}{t-1}}
\end{eqnarray*}
by the compositional formula for exponential generating functions.
$\Box$ \medskip

\begin{theorem}\label{th:Tn}
The threshold arrangement ${\mathcal T}_n$ in $\RR^n$ consists of
the hyperplanes $x_i + x_j = 0$, for $1 \leq i < j \leq n$. For
all $n \geq 0$ we have
$$
\Chi_{{\mathcal T}_n} (q,t) = \sum_G q^{bc(G)} (t-1)^{e(G)},
$$
where the sum is over all graphs $G$ on $[n]$. Here $bc(G)$ is
the number of connected components of $G$ which are bipartite,
and $e(G)$ is the number of edges of $G$. Also,
$$
\sum_{n \geq 0} \Chi_{{\mathcal T}_n}(q,t) \frac{x^n}{n!} =
\left(\sum_{n \geq 0} \sum_{k=0}^n {n \choose k}
t^{k(n-k)}\frac{x^n}{n!} \right)^{\frac{q-1}2} \left(\sum_{n \geq
0} t^{n \choose 2} \frac{x^n}{n!} \right).
$$
\end{theorem}

\noindent \emph{Proof.} Once again, the proof of the first claim
is easier using the definition of the coboundary polynomial.
Every subarrangement $\B$ of $\T_n$ is central, and we can assign
to it a graph $G_{\B}$ as in the proof of Theorem \ref{th:An'}. In
view of (\ref{eq:cob}), we only need to check that $r(\B) = n -
bc(G_{\B})$ and $|\B| = e(G_{\B})$. The second claim is trivial.
To prove the first one, we show that $\dim(\cap \B)=bc(G_{\B})$.

Consider a point $p$ in $\cap \B$. We know that, if $ab$ is an
edge in $G_{\B}$, then $p_a = - p_b$. If vertex $i$ is in a
connected component $C$ of $G_{\B}$, then the value of $p_i$
determines the value of $p_j$ for all $j$ in $C$: $p_j = p_i$ if
there is a path of even length between $i$ and $j$, and $p_j =
-p_i$ if there is a path of odd length between $i$ and $j$. If $C$
is bipartite, this determines the values of the $p_j$'s
consistently. If $C$ is not bipartite, take a cycle of odd length
and a vertex $k$ in it. We get that $p_k = -p_k$, so $p_k=0$;
therefore we must have $p_j=0$ for all $j \in C$.

Therefore, to specify a point $p$ in $\cap \B$, we split $G_{\B}$
into its connected components. We know that $p_i = 0$ for all $i$
in connected components which are not bipartite. To determine the
remaining coordinates of $p$ we have to specify the value of
$p_j$ for exactly one $j$ in each bipartite connected component.
Therefore $\dim(\cap \B)=bc(G_{\B})$, as desired.

From this point, it is possible to prove the second claim of
Theorem \ref{th:Tn} using the compositional formula for
exponential generating functions, in the same way that we proved
Theorem \ref{th:An'}. However, the work involved is considerable,
and it is much simpler to use our finite field method, Theorem
\ref{th:finfield}, in this case. The proof that we obtain is very
similar to the proofs of Theorems \ref{th:An}, \ref{th:Bn} and
\ref{th:Dn}, so we omit the details. $\Box$

\subsection{Deformations of the braid arrangement.}

A \emph{deformation of the braid arrangement} is an arrangement in
$\RR^n$ consisting of the hyperplanes $x_i - x_j = a_{ij}^{(1)},
\ldots, a_{ij}^{(k_{ij})}$ for $1 \leq i < j \leq n$, where the
$k_{ij}$ are non-negative integers, and the $a_{ij}^{(r)}$ are
real numbers. Such arrangements have been studied extensively by
Athanasiadis \cite{At98} and Postnikov and Stanley \cite{Po00}.
In this section we study their coboundary polynomials.

\begin{theorem}\label{th:esa}
Let $\E = (\E_0, \E_1, \ldots)$ be a sequence of arrangements
satisfying the following properties: \footnote{Such a sequence is
called an \emph{exponential sequence of arrangements.}}
\begin{enumerate}
\item
$\E_n$ is an arrangement in $\kk^n$, for some fixed field $\kk$.
\item
Every hyperplane in $\E_n$ is parallel to some hyperplane in the
braid arrangement $\A_n$.
\item
For any subset $S$ of $[n]$, the subarrangement $\E_n^S \subseteq
\E_n$, which consists of the hyperplanes in $\E_n$ of the form
$x_i-x_j=c$ with $i,j \in S$, is isomorphic to the arrangement
$\E_{|S|}$.
\end{enumerate}
Then
$$
1 + q \sum_{n \geq 1} \Chi_{\E_n}(q,t) \frac{x^n}{n!} = \left(\sum_{n
\geq 0} \Chi_{\E_n}(1,t) \frac{x^n}{n!}\right)^q.
$$
\end{theorem}

\medskip

The special case $t=0$ of this result is due to Stanley
\cite[Theorem 1.2]{St96}; we omit the proof, which is an easy
extension of his.

\medskip

The most natural examples of exponential sequences of
arrangements are the following. Fix a set $A$ of $k$ distinct
integers $a_1< \ldots< a_k$. Let $\E_n$ be the arrangement in
$\RR^n$ consisting of the hyperplanes
\begin{equation}\label{eq:En}
x_i - x_j = a_1, \ldots, a_k \qquad 1 \leq i < j \leq n.
\end{equation}
Then $(\E_0, \E_1, \ldots)$ is an exponential sequence of
arrangements and Theorem \ref{th:esa} applies to this case. In
fact, we can say much more about this type of arrangement.

After proving the results in this section, we found out that
Postnikov and Stanley \cite{Po00} had used similar techniques in
computing the characteristic polynomials of these types of
arrangements. Therefore, for consistency, we will use the
terminology that they introduced.

\begin{defin}
A \emph{graded graph} is a triple $G = (V_G, E_G, h_G)$, where
$V_G$ is a linearly ordered set of vertices (usually $V_G=[n]$),
$E_G$ is a set of (non-oriented) edges, and $h_G$ is a function
$h_G:V \rightarrow \NN$, called a \emph{grading}.
\end{defin}

We will drop the subscripts when the underlying graded graph is
clear. We will refer to $h(v)$ as the \emph{height} of vertex
$v$. The \emph{height} of $G$, denoted $h(G)$, is the largest
height of a vertex of $G$.

\begin{defin}
Let $G$ be a graded graph and $r$ be a  non-negative integer. The
\emph{$r$-th level} of $G$ is the set of vertices $v$ such that
$h(v)=r$. $G$ is \emph{planted} if each one of its connected
components has a vertex on the $0$-th level.
\end{defin}

\begin{defin}
If $u<v$ are connected by edge $e$ in a graded graph $G$, the
\emph{slope} of $e$ is $s(e) = h(u)-h(v)$. $G$ is an
\emph{$A$-graph} if the slopes of all edges of $G$ are in
$A=\{a_1, \ldots, a_k\}$.
\end{defin}

Recall that, for a graph $G$, we let $e(G)$ be the number of
edges and $c(G)$ be the number of connected components of $G$. We
also let $v(G)$ be the number of vertices of $G$.

\begin{proposition}\label{pr:En}
Let $\E_n$ be the arrangement (\ref{eq:En}). Then, for $n \geq 1$,
$$
q \, \Chi_{\E_n}(q,t) = \sum_G q^{c(G)}(t-1)^{e(G)},
$$
where the sum is over all planted graded $A$-graphs on $[n]$.
\end{proposition}

\noindent \emph{Proof.} We associate to each planted graded
$A$-graph $G = (V, E, h)$ on $[n]$ a central subarrangement
$\A_G$ of $\E_n$. It consists of the hyperplanes $x_i-x_j =
h(i)-h(j)$, for each $i < j$ such that $ij$ is an edge in $G$.
This is a subarrangement of $\E_n$ because $h(i)-h(j)$, the slope
of edge $ij$, is in $A$. It is central because the point $(h(1),
\ldots, h(n)) \in \RR^n$ belongs to all these hyperplanes.

This is in fact a bijection between planted graded $A$-graphs on
$[n]$ and central subarrangements of $\E_n$. To see this, take a
central subarrangement $\A$. We will recover the planted graded
$A$-graph $G$ that it came from. For each pair $(i,j)$ with $1
\leq i < j \leq n$, $\A$ can have at most one hyperplane of the
form $x_i-x_j=a_t$. If this hyperplane is in $\A$, we must put
edge $ij$ in $G$, and demand that the heights $h(i)$ and $h(j)$
satisfy $h(i)-h(j)=a_t$. When we do this for all the hyperplanes
in $\A$, the height requirements that we introduce are consistent,
because $\A$ is central. However, these requirements do not fully
determine the heights of the vertices; they only determine the
relative heights within each connected component of $G$. Since we
want $G$ to be planted, we demand that the vertices with the
lowest height in each connected component of $G$ should have
height $0$. This does determine $G$ completely, and clearly $\A =
\A_G$.

\medskip

\begin{figure}[tbp]
\centering
\includegraphics[width=2.5in]{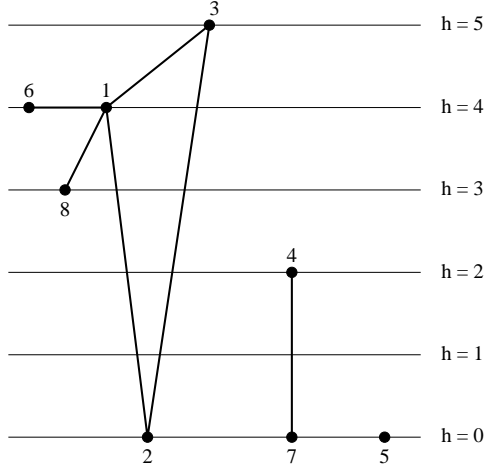}
\caption{The planted graded $A$-graph corresponding to a
subarrangement of $\E_8$.} \label{fig:Agraph}
\end{figure}

\noindent \emph{Example.} Consider an arrangement $\E_8$ in
$\RR^8$, with a subarrangement  consisting of the hyperplanes
$x_1-x_2=4, x_1-x_3=-1, x_1-x_6=0, x_1-x_8=1, x_2-x_3=-5$ and
$x_4-x_7=2$. Figure \ref{fig:Agraph} shows the planted graded
$A$-graph corresponding to this subarrangement.

\medskip

With this bijection in hand, and keeping (\ref{eq:cob}) in mind,
it remains to show that $r(\A_G) = n - c(G)$ and $|\A_G|=e(G)$.
The second of these claims is trivial. We omit the proof of the
first one which is very similar to, and simpler than, that of
$r(\B) = n - bc(G_{\B})$ in our proof of Theorem \ref{th:Tn}.
$\Box$

\medskip

\begin{theorem}\label{th:deform} Let $\E_n$ be the
arrangement (\ref{eq:En}), and let
\begin{equation}\label{eq:defAr}
A_r(t,x) = \sum_{n \geq 0} \left(\sum_{f:[n]\rightarrow [r]}
t^{a(f)}\right) \frac{x^n}{n!},
\end{equation}
where $a(f)$ denotes the number of pairs $(i,j)$ with $1 \leq i <
j \leq n$ such that $f(i) - f(j) \in A$. Then
\begin{equation}\label{eq:lim}
1 + q \sum_{n \geq 1} \Chi_{\E_n}(q,t) \frac{x^n}{n!} = \left(\lim_{r
\rightarrow \infty} \frac{A_r(t,x)}{A_{r-1}(t,x)}\right)^q.
\end{equation}
\end{theorem}

\medskip

\noindent \emph{Remark.} The limit in (\ref{eq:lim}) is a
limit in the sense of convergence of formal power series. Let $F_1(t,x),
F_2(t,x), \ldots$ be a sequence of formal power series. We say
that $\lim_{n \rightarrow \infty} F_n(t,x) = F(t,x)$ if, for all
$a$ and $b$, the coefficient of $t^ax^b$ in $F_n(t,x)$ is equal to
the coefficient of $t^ax^b$ in $F(t,x)$ for all $n$ larger than
some constant $N(a,b)$. For more information on this notion of
convergence, see \cite[Section 1.1]{St86} or
\cite{Ni69}.

\medskip

\noindent \emph{Proof of Theorem \ref{th:deform}.} First we prove
that
\begin{equation}\label{eq:Ar}
A_r(t,x) = \sum_G (t-1)^{e(G)} \frac{x^{v(G)}}{v(G)!}
\end{equation}
where the sum is over all graded $A$-graphs $G$ of height less
than $r$. The coefficient of $\frac{x^n}{n!}$ in the right-hand
side of (\ref{eq:Ar}) is $\sum_G (t-1)^{e(G)}$, summing
over all graded $A$-graphs $G$ on $[n]$ with height less than $r$.
We have
\begin{eqnarray*}
\sum_G (t-1)^{e(G)}&=& \sum_{h:[n] \rightarrow [0, r-1]} \,\,
\sum_{\stackrel{G \, \mathrm{such \, that}} {h_G = h}} (t-1)^{e(G)} \\
&=& \sum_{h:[n] \rightarrow [0, r-1]} (1 + (t-1))^{a(h)} \\
&=& \sum_{f:[n] \rightarrow [r]} t^{a(f)}
\end{eqnarray*}
The only tricky step here is the second: if we want all graded
$A$-graphs $G$ on $[n]$ with a specified grading $h$, we need to
consider the possible choices of edges of the graph. Any edge
$ij$ can belong to the graph, as long as $h(i)-h(j) \in A$, so
there are $a(h)$ possible edges.

Equation (\ref{eq:Ar}) suggests the following definitions. Let
$$
B_r(t,x) = \sum_G t^{e(G)} \frac{x^{v(G)}}{v(G)!}
$$
where the sum is over all \emph{planted} graded $A$-graphs $G$ of
height less than $r$, and let
$$
B(t,x) = \sum_G t^{e(G)} \frac{x^{v(G)}}{v(G)!}
$$
where the sum is over all \emph{planted} graded $A$-graphs $G$.

The equation
\begin{equation}\label{eq:EB}
1 + q \sum_{n \geq 1} \Chi_{\E_n}(q,t) \frac{x^n}{n!} = B(t-1,x)^q,
\end{equation}
follows from Proposition \ref{pr:En}, using either Theorem
\ref{th:esa} or the compositional formula for exponential
generating functions.

Now we claim that $B(t,x) = \lim_{r \rightarrow \infty}
B_r(t,x)$. Notice that, in a planted graded $A$-graph $G$ with $e$
edges and $v$ vertices, each vertex has a path of length at most
$v$ which connects it to a vertex on the $0$-th
level. Recalling that $a_1<\ldots<a_k$ we see that $h(G) \leq v
\cdot \max(-a_1, a_k)$, so the coefficients of $t^e\frac{x^v}{v!}$
in $B_r(t,x)$ and $B(t,x)$ are equal for $r > v \cdot
\max(-a_1, a_k)$.

With a little bit of care, it then follows easily that
\begin{equation}\label{eq:limit}
B(t-1,x) = \lim_{r \rightarrow \infty} B_r(t-1,x).
\end{equation}
Here it is necessary to check that $B(t-1,x)$ is, indeed, a formal
power series. This follows from the observation that the
coefficient of $\frac{x^n}{n!}$ in $B(t,x)$ is a polynomial in
$t$ of degree at most ${n \choose 2}$. We know that for some
formal power series $f(t)$ (like $e^t$, for example), $f(t-1)$ is
not a well-defined formal power series.  In our case, however,
this is not a problem and (\ref{eq:limit}) is valid. Once again,
see \cite[Section 1.1]{St86} for more information on these
technical details.

Next, we show that
\begin{equation}\label{eq:AB}
B_r(t-1,x) = A_r(t,x)/A_{r-1}(t,x)
\end{equation}
or, equivalently, that $A_r(t,x) = B_r(t-1,x) A_{r-1}(t,x).$ The
multiplication formula for exponential generating functions
(\cite[Proposition 5.1.1]{St99}) and (\ref{eq:Ar}) give us a
combinatorial interpretation of this identity. We need to show
that the ways of putting the structure of a graded $A$-graph $G$
with $h(G)<r$ on $[n]$ can be put in correspondence with the ways
of doing the following: first splitting $[n]$ into two disjoint
sets $S_1$ and $S_2$, then putting the structure of a
\emph{planted} graded $A$-graph $G_1$ with $h(G_1)<r$ on $S_1$,
and then putting the structure of a graded $A$-graph $G_2$ with
$h(G_2)<r-1$ on $S_2$. We also need that, in that correspondence,
$(t-1)^{e(G)} = (t-1)^{e(G_1)}(t-1)^{e(G_2)}$.

\begin{figure}[tbp]
\centering
\includegraphics[width=5in]{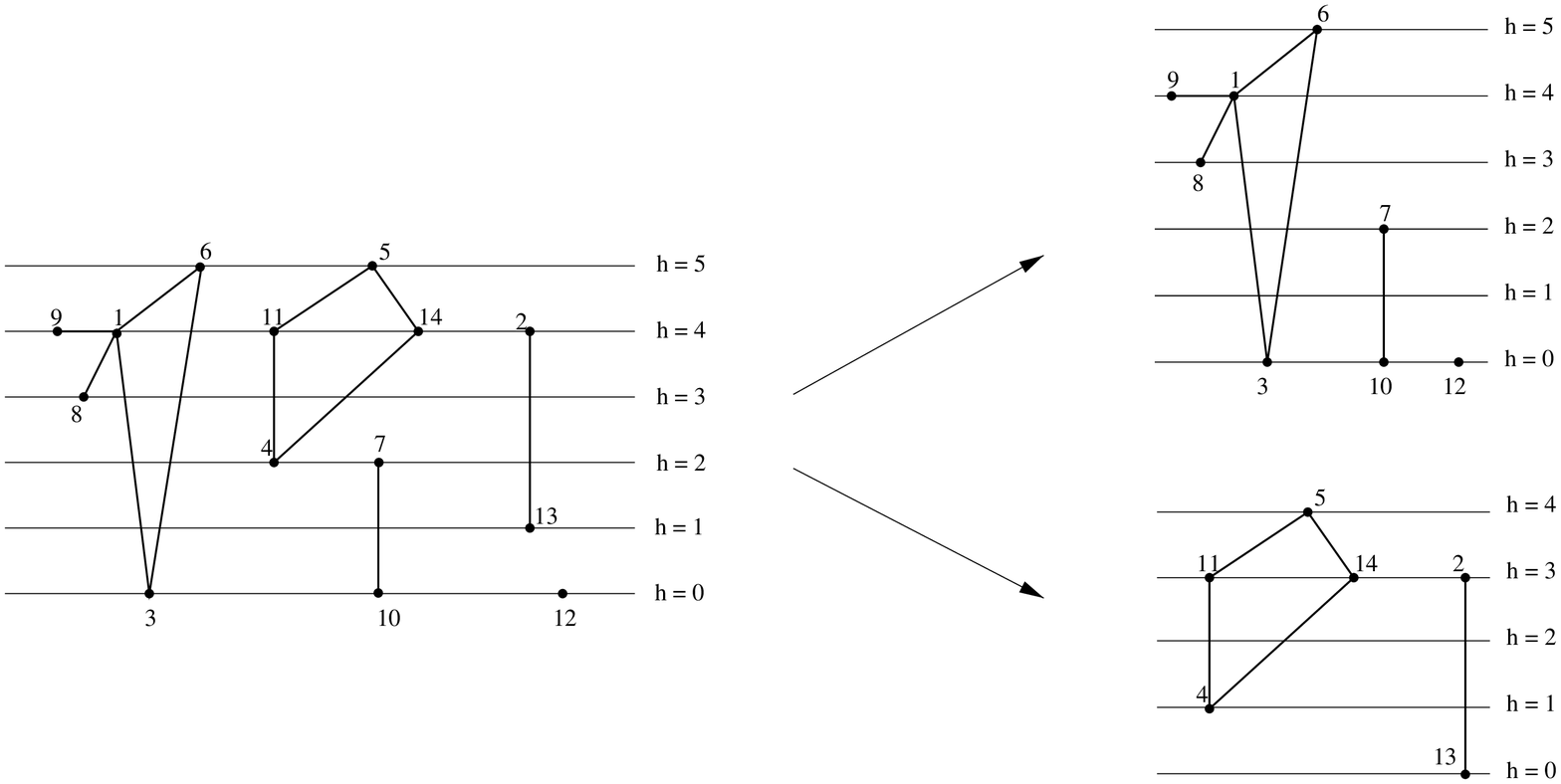}
\caption{The decomposition of a graded $A$-graph.}
\label{fig:decomp}
\end{figure}

We do this as follows. Let $G$ be a graded $A$-graph $G$ with
$h(G)<r$. Let $G_1$ be the union of the connected components of
$G$ which contain a vertex on the $0$-th level. Put a grading on
$G_1$ by defining $h_{G_1}(v) = h_G(v)$ for $v \in G_1$. Let
$G_2 = G - G_1$. It is clear that $h_G(v) \geq 1$ for all $v \in
G_2$; therefore we can put a grading on $G_2$ by defining
$h_{G_2}(v) = h_G(v) - 1$ for $v \in G_2$. $G_1$ is a planted graded
$A$-graph with $h(G_1)<r$, and $G_2$ is a graded $A$-graph with
$h(G_2)<r-1$. Figure \ref{fig:decomp} illustrates this
decomposition with an example.

Our map from $G$ to a pair $(G_1, G_2)$ is a one-to-one
correspondence. Any pair $(G_1,G_2)$, with $G_1$ planted of
height less than $r$ and $G_2$ of height less than $r-1$, arises
from a decomposition of some $G$ of height less than $r$ in this
way. It is clear how to recover $G$ from $G_1$ and $G_2$. Also,
it is clear from the construction of the correspondence that
$(t-1)^{e(G)} = (t-1)^{e(G_1)}(t-1)^{e(G_2)}$. This completes the
proof of (\ref{eq:AB}).

Now we just have to put together (\ref{eq:EB}), (\ref{eq:limit})
and (\ref{eq:AB}) to complete the proof of Theorem
\ref{th:deform}. $\Box$

\medskip

The \emph{Catalan arrangement} $C_n$ in $\RR^n$ consists of the
hyperplanes
$$
x_i - x_j = -1, 0, 1 \qquad 1 \leq i < j \leq n.
$$
When the arrangement in Theorem \ref{th:deform} is a
subarrangement of the Catalan arrangement, we can say more about
the power series $A_r$ of (\ref{eq:defAr}). Let
$$
A(t,x,y) = \sum_r A_r(t,x)y^r = \sum_{n \geq 0} \sum_{r \geq 0}
\left(\sum_{f:[n]\rightarrow [r]} t^{a(f)}\right) \frac{x^n}{n!}
y^r
$$
and let
\begin{equation}\label{eq:defS}
S(t,x,y) = \sum_{n \geq 0} \sum_{r \geq 0}
\left(\sum_{f:[n]\rightarrow [r]} t^{a(f)}\right) \frac{x^n}{n!}
y^r
\end{equation}
where the inner sum is over all \emph{surjective} functions $f:[n]
\rightarrow [r]$. The following proposition reduces the computation
of $A(t,x,y)$ to the computation of $S(t,x,y)$, which is often
easier in practice.

\begin{proposition}\label{pr:AS}
If $A \subseteq \{-1,0,1\}$ in the notation of Theorem
\ref{th:deform}, we have
$$
A(t,x,y) = \frac{S(t,x,y)}{1 - yS(t,x,y)}
$$
\end{proposition}

\noindent \emph{Proof.} Once again, we think of this as an
identity about exponential generating functions in the variable
$x$. Fix $n,r$, and $f:[n] \rightarrow [r]$. Let the image of $f$
be $\{1 , \ldots , m_1-1\} \cup \{m_1+1 , \ldots , m_1+m_2-1\} \cup
\cdots \cup \{m_1 + \cdots + m_{k-1}+1 , \ldots , m_1 + \cdots +
m_k-1\} = M_1 \cup \cdots \cup M_k$, so that $[r] - \Im f = \{m_1,
m_1+m_2, \ldots, m_1 + \cdots + m_{k-1}\}$. Here $m_1, \ldots,
m_k$ are arbitrary positive integers such that $m_1 + \cdots +
m_k - 1 = r$. For $1 \leq i \leq k$, let $f_i$ be the restriction
of $f$ to $f^{-1}(M_i)$; it maps $f^{-1}(M_i)$ surjectively to
$M_i$. Then we can ``decompose" $f$ in a unique way
into the $k$ \emph{surjective} functions $f_1, \ldots, f_k$. The
weight $w(f)$ corresponding to $f$ in $A(t,x,y)$ is $t^{a(f)}
y^r$, while the weight $w(f_i)$ corresponding to $f_i$ in
$S(t,x,y)$ is $t^{a(f_i)} y^{m_i - 1}$.

Now observe that $a(f) = a(f_1) + \cdots + a(f_k)$: whenever we
have a pair of numbers $1 \leq i < j \leq n$ counted by $a(f)$,
since $f(i) - f(j) \in \{-1,0,1\}$, we know that $f(i)$ and
$f(j)$ must be in the same $M_h$. Therefore $i$ and $j$ are in
the same $f^{-1}(M_h)$, and this pair is also counted in $a(f_h)$.
We also have that $r = (m_1 - 1) + \cdots + (m_k - 1) + (k-1)$.
Therefore $w(f) = w(f_1) \cdots w(f_k) y^{k-1}$. It follows from the
compositional formula for exponential generating functions that
\begin{eqnarray*}
A(t,x,y) &=& \sum_{k \geq 1} S(t,x,y)^k y^{k-1} \\
&=& \frac{S(t,x,y)}{1 - yS(t,x,y)}
\end{eqnarray*}
as desired. $\Box$

\medskip

Considering the different subsets of $\{-1,0,1\}$, we get six
non-isomorphic subarrangements of the Catalan arrangement. They
come from the subsets $\emptyset$, $\{0\}$, $\{1\}$, $\{0,1\}$,
$\{-1,1\}$ and $\{-1,0,1\}$. The corresponding subarrangements are
the empty arrangement, the braid arrangement, the Linial
arrangement, the Shi arrangement, the semiorder arrangement, and
the Catalan arrangement, respectively. The empty arrangement is
trivial, and the braid arrangement was already treated in detail
in Section \ref{sec:cox}. We now have a technique that lets us
talk about the remaining four arrangements under the same
framework. We will do this in the remainder of this chapter.

\subsubsection{The Linial arrangement}

The Linial arrangement $\L_n$ consists of the hyperplanes $x_i -
x_j = 1$ for $1 \leq i < j \leq n$. This arrangement was first
considered by Linial and Ravid. It was later studied
by Athanasiadis \cite{At96} and Postnikov and Stanley
\cite{Po00}, who independently computed the characteristic
polynomial of $\L_n$:
$$
\chi_{\L_n}(q) = \frac q{2^n} \sum_{k=0}^n {n \choose k}
(q-k)^{n-1}.
$$
They also put the regions of $\L_n$ in bijection with several
different sets of combinatorial objects. Perhaps the simplest
such set is the set of \emph{alternating trees} on $[n+1]$: the
trees such that every vertex is either larger or smaller
than all its neighbors.

Now we present the consequences of Proposition \ref{pr:En},
Theorem \ref{th:deform} and Proposition \ref{pr:AS} for the Linial
arrangement. Recall that a poset $P$ on $[n]$ is \emph{naturally
labeled} if $i<j$ in $P$ implies $i<j$ in $\ZZ^+$.

\begin{proposition}
For all $n \geq 1$ we have
$$
q \, \Chi_{\L_n} (q,t) = \sum_P q^{c(P)} (t-1)^{e(P)}
$$
where the sum is over all naturally labeled, graded posets $P$ on
$[n]$. Here $c(P)$ and $e(P)$ denote the number of components and
edges of the Hasse diagram of $P$, respectively.
\end{proposition}

\noindent \emph{Proof.} There is an obvious bijection between
Hasse diagrams of naturally labeled graded posets on $[n]$ and
planted graded $\{1\}$-graphs on $[n]$. The result then follows
immediately from Proposition \ref{pr:En}. $\Box$ \medskip

\begin{theorem} \label{th:Ln}
Let
$$
A_r(t,x) = \sum_{n \geq 0} \left(\sum_{f:[n]\rightarrow [r]}
t^{id(f)}\right) \frac{x^n}{n!}.
$$
where $id(f)$ denotes the number of \emph{inverse descents} of
the word $f(1) \ldots f(n)$: the number of pairs $(i,j)$ with $1
\leq i < j \leq n$ such that $f(i) - f(j) = 1$. Then
$$
1 + q \sum_{n \geq 1} \Chi_{\L_n}(q,t) \frac{x^n}{n!} = \left(
\lim_{r \rightarrow \infty} \frac{A_r(t,x)}{A_{r-1}(t,x)}\right)^q.
$$
\end{theorem}

\noindent \emph{Proof.} This is immediate from Theorem
\ref{th:deform}. $\Box$

\medskip

Recall that the \emph{descents} of a permutation $\sigma =
\sigma_1 \ldots \sigma_r \in S_r$ are the indices $i$ such that
$\sigma_i > \sigma_{i+1}$. For more information about descents,
see for example \cite[Section 1.3]{St86}.

We call $id(f)$ the number of inverse descents, because they
generalize descents in the following way. If $\pi:[r]\rightarrow
[r]$ is a permutation, then $id(\pi)$ is the number of descents
of the permutation $\pi^{-1}$. If, similarly, we consider the list
of sets $f^{-1}(1), \ldots, f^{-1}(r)$, then $id(f)$ counts the
number of occurrences of an $x \in f^{-1}(i)$ and a $y \in
f^{-1}(i+1)$ such that $x>y$.

It would be nice to compute the polynomials $A_r(t,x)$ above
explicitly. We have not been able to do this. However, the special
case $t=0$ is of interest; recall that the characteristic
polynomial of $\L_n$ is $\chi_{\L_n}(q) = q\Chi_{\L_n}(q,0)$. In
that case, we obtain the following result.

\begin{theorem}\label{th:charLn}
Let
\begin{equation}\label{eq:ALn}
\frac{1+ye^{x(1+y)}}{1-y^2e^{x(1+y)}} = \sum_{r \geq 0} A_r(x)y^r.
\end{equation}
Then we have
$$
\sum_{n \geq 0} \chi_{\L_n}(q) \frac{x^n}{n!} = \left(\lim_{r
\rightarrow \infty} \frac{A_r(x)}{A_{r-1}(x)}\right)^q.
$$
In particular, if $f_n$ is the number of alternating trees on
$[n+1]$, we have
$$
\sum_{n \geq 0} (-1)^n f_n \frac{x^n}{n!} = \lim_{r
\rightarrow \infty} \frac{A_{r-1}(x)}{A_r(x)}.
$$
\end{theorem}

\noindent \emph{Proof.} In view of Theorem \ref{th:Ln} and
Proposition \ref{pr:AS}, we compute $S(0,x,y)$. From
(\ref{eq:defS}), the coefficient of $\frac{x^n}{n!} y^r$ in
$S(0,x,y)$ is equal to the number of surjective functions $f:[n]
\rightarrow [r]$ with no inverse descents. These are just the
non-decreasing surjective functions $f:[n] \rightarrow [r]$. For
$n \geq 1$ and $r \geq 1$ there are ${n-1 \choose r-1}$ such
functions, and for $n=r=0$ there is one such function. In the
other cases there are none. Therefore
\begin{eqnarray*}
S(0,x,y) & = & 1 + \sum_{n \geq 1} \sum_{r \geq 1} {n-1 \choose r-1}
\frac{x^n}{n!} y^r \\
& = & 1 + \sum_{n \geq 1} \frac{x^n}{n!} y(1+y)^{n-1} \\
& = & \frac{1 +ye^{x(1+y)}}{1+y}.
\end{eqnarray*}
Proposition \ref{pr:AS} then implies that
$$
A(0,x,y) = \frac{1+ye^{x(1+y)}}{1-y^2e^{x(1+y)}},
$$
in agreement with (\ref{eq:ALn}), and the theorem follows. $\Box$

\subsubsection{The Shi arrangement}

The Shi arrangement $\S_n$ consists of the hyperplanes $x_i - x_j
= 0, 1$ for $1 \leq i < j \leq n$. Shi (\cite[Chapter 7]{Sh86},
\cite{Sh87}) first considered this arrangement, and showed that
it has $(n+1)^{n-1}$ regions. Headley (\cite[Chapter VI]{He94},
\cite{He97}) later computed the characteristic polynomial of
$\S_n$:
$$
\chi_{\S_n}(q) = q(q-n)^{n-1}.
$$
Stanley \cite{St96},\cite{St98} gave a nice bijection between
regions of the Shi arrangement and parking functions of length
$n$. Parking functions were first introduced by Konheim and Weiss
\cite{Ko66}; for more information about them, see \cite[Exercise 5.49]
{St99}.

For the Shi arrangement, we can say the following.

\begin{theorem}\label{th:charSn}
Let
$$
A_r(x) = \sum_{n=0}^r (r-n)^n \frac{x^n}{n!}.
$$
Then we have
$$
\sum_{n \geq 0} \chi_{\S_n}(q) \frac{x^n}{n!} = \left(\lim_{r
\rightarrow \infty} \frac{A_r(x)}{A_{r-1}(x)}\right)^q.
$$
In particular, we have
$$
\sum_{n \geq 0} (-1)^n (n+1)^{n-1} \frac{x^n}{n!} = \lim_{r
\rightarrow \infty} \frac{A_{r-1}(x)}{A_r(x)}.
$$
\end{theorem}

\noindent \emph{Proof.} We proceed in the same way that we did in
Theorem \ref{th:charLn}. In this case, we need to compute the
number of surjective functions $f:[n] \rightarrow [r]$ such that
$f(i)-f(j)$ is never equal to $0$ or $1$ for $i<j$. These are
just the surjective, strictly increasing functions. There is only
one of them when $n=r$, and there are none when $n \neq r$.
Hence
$$
S(0,x,y) = \sum_{n \geq 0} \frac{x^n}{n!} y^n = e^{xy}.
$$
The rest follows easily by computing $A(0,x,y)$ and $A_r(x)$
explicitly. $\Box$

\subsubsection{The semiorder arrangement}

The semiorder arrangement $\I_n$ consists of the hyperplanes $x_i
- x_j = -1, 1$ for $1 \leq i < j \leq n$. A \emph{semiorder on
$[n]$} is a poset $P$ on $[n]$ for which there exist $n$ unit
intervals $I_1, \ldots, I_n$ of $\RR$ , such that $i<j$ in $P$ if
and only if $I_i$ is disjoint from $I_j$ and to the left of it. It
is known \cite{Sc58} that a poset is a semiorder if and only if it
does not contain a subposet isomorphic to $\bf{3}+\bf{1}$ or
$\bf{2}+\bf{2}$. We are interested in semiorders because the
number of regions of $\I_n$ is equal to the number of semiorders
on $[n]$, as shown in \cite{Po00} and \cite{St96}.

\begin{theorem}\label{charIn}
Let
$$
\frac{1-y+ye^x}{1-y+y^2-y^2e^x}= \sum_{r \geq 0} A_r(x)y^r.
$$
Then we have
$$
\sum_{n \geq 0} \chi_{\I_n}(q) \frac{x^n}{n!} = \left( \lim_{r
\rightarrow \infty} \frac{A_r(x)}{A_{r-1}(x)}\right)^q.
$$
In particular, if $i_n$ is the number of semiorders on $[n]$, we
have
$$
\sum_{n \geq 0} (-1)^n i_n \frac{x^n}{n!} = \lim_{r \rightarrow
\infty} \frac{A_{r-1}(x)}{A_r(x)}.
$$
\end{theorem}

\noindent \emph{Proof.} In this case, $S(0,x,y)$ counts
surjective functions $f:[n] \rightarrow [r]$ such that
$f(i)-f(j)$ is never equal to $1$ for $i \neq j$. Such a function
has to be constant; so it can only exist (and is unique) if $n
\geq 1$ and $r=1$ or if $n=r=0$. Thus
$$
S(0,x,y) = 1 + (e^x-1)y
$$
and the rest follows easily. $\Box$

\subsubsection{The Catalan arrangement}

The Catalan arrangement $C_n$ consists of the hyperplanes $x_i -
x_j = -1,0,1$ for $1 \leq i < j \leq n$. Stanley \cite{St96}
observed that the number of regions of this arrangement is $n!
C_n$, where $C_n= \frac1{n+1} {2n \choose n}$ is the $n$-th Catalan
number. For (much) more information on the Catalan numbers, see
\cite[Chapter 6]{St99}, especially Exercise 6.19.

\begin{theorem}\label{charCn}
Let
$$
A_r(x) = \sum_{n = 0}^{\lfloor \frac{r+1}2 \rfloor} {r-n+1 \choose n}x^n.
$$
Then we have
$$
\sum_{n \geq 0} \chi_{C_n}(q) \frac{x^n}{n!} = \left( \lim_{r
\rightarrow \infty} \frac{A_r(x)}{A_{r-1}(x)}\right)^q.
$$
In particular,
\begin{equation}
\frac{\sqrt{1+4x}-1}{2x} = \sum_{n \geq 0} (-1)^n C_n x^n =
\lim_{r \rightarrow \infty} \frac{A_{r-1}(x)}{A_r(x)}.
\label{eq:fibonacci}
\end{equation}
\end{theorem}

\noindent \emph{Proof.} There are no surjective functions $f:[n]
\rightarrow [r]$ such that $f(i)-f(j)$ is never equal to $-1,0$
or $1$ for $i \neq j$, unless $n=r=0$ or $n=r=1$. Thus $S(x,y,0)
= 1+xy$. The rest of the proof is easy. $\Box$

\medskip

The polynomial $A_r(x)$ is a simple transformation of the
\emph{Fibonacci polynomial}. The number of words of length $r$,
consisting of $0$'s and $1$'s, which do not contain two
consecutive $1$'s, is equal to $F_{r+2}$, the $(r+2)$-th Fibonacci
number. It is easy to see that the polynomial $A_r(x)$ counts
those words according to the number of $1$'s they contain. In
particular, $A_r(1)=F_{r+2}$.

We close with an amusing observation. Irresponsibly plugging $x=1$
into (\ref{eq:fibonacci}) \footnote{We are not justified in doing
this!}, we obtain an unconventional ``proof" of the rate of growth
of Fibonacci numbers:
$$
\frac{\sqrt{5}-1}2 = \lim_{r \rightarrow \infty}
\frac{F_{r-1}}{F_r}.
$$

\section{Acknowledgments.}

The present work is Chapter 2 of the author's Ph.D. thesis
\cite{Ar02}. I would like to thank my advisor, Richard Stanley,
for introducing me to the topic of hyperplane arrangements, and
for asking some of the questions which led to this investigation.
I am also grateful to Ira Gessel and Vic Reiner for very helpful
discussions on this subject.

\end{document}